\documentclass[12pt]{amsart}
\usepackage[utf8]{inputenc}

\title{A note on uncountably chromatic graphs}
\author[N.~Bowler]{Nathan Bowler}
\author[M.~Pitz]{Max Pitz}
\address{Universit\"at Hamburg, Department of Mathematics, Bundesstrasse 55 (Geomatikum), 20146 Hamburg, Germany}
\email{\{nathan.bowler, max.pitz\}@uni-hamburg.de}

\date{}



\usepackage{amsmath,amssymb,amsthm}  
\usepackage{mathtools}
\usepackage{tikz}
\usepackage{comment}
\usepackage{enumitem}

\usepackage{xcolor} 	
\usepackage[unicode]{hyperref}
\hypersetup{
	colorlinks,
    linkcolor={red!60!black},
    citecolor={green!60!black},
    urlcolor={blue!60!black},
}
\usepackage[abbrev, msc-links]{amsrefs} 

\usepackage[utf8]{inputenc}
\usepackage[T1]{fontenc}
\usepackage{lmodern}
\usepackage[babel]{microtype}
\usepackage[english]{babel}

\usepackage{geometry}
\usepackage{enumitem}

\let\polishlcross=\l
\def\l{\ifmmode\ell\else\polishlcross\fi}

\let\theta=\vartheta
\let\rho=\varrho
\let\phi=\varphi

\def\NN{\mathbb N}

\newcommand{\parentheses}[1]{{\left( {#1} \right)}}

\newcommand{\p}{\parentheses}

\newcommand{\Set}[1]{{\left\lbrace {#1} \right\rbrace}}

\def\set#1:#2{\Set{{#1} \colon {#2}}}

\DeclareFontFamily{U}  {MnSymbolC}{}
\DeclareSymbolFont{MnSyC}         {U}  {MnSymbolC}{m}{n}
\DeclareFontShape{U}{MnSymbolC}{m}{n}{
    <-6>  MnSymbolC5
   <6-7>  MnSymbolC6
   <7-8>  MnSymbolC7
   <8-9>  MnSymbolC8
   <9-10> MnSymbolC9
  <10-12> MnSymbolC10
  <12->   MnSymbolC12}{}
\DeclareMathSymbol{\powerset}{\mathord}{MnSyC}{180}

\theoremstyle{plain}
\newtheorem*{thm}{Theorem}
\newtheorem*{clm}{Claim}
\begin{document}
\begin{abstract}
We present an elementary construction of an uncountably chromatic graph without uncountable, infinitely connected subgraphs. 
\end{abstract}

\maketitle

\section{Introduction}

Erd\H{o}s and Hajnal asked in 1985 whether every graph of uncountable chromatic number has an infinitely connected, uncountably chromatic subgraph. In 1988 and 2013, P.~Komj\'ath gave  consistent negative answers \cite{komjath1988consistency,komjath2013note}: He first constructed an uncountably chromatic graph without infinitely connected, uncountably chromatic subgraphs; and later an uncountably chromatic graph without any uncountable, infinitely connected subgraph.
In 2015, D.~Soukup gave the first ZFC construction of 
such a graph \cite{soukup2015trees}.
Soukup even produces an uncountably chromatic graph $G$ in which every uncountable set of vertices contains two points that are connected by only finitely many independent paths in $G$.
In this note we present a short, elementary example for Soukup's result.

\section{The example}
Let $\NN = \{1,2,3,\ldots\}$. 
For a countable ordinal $\alpha$, write $T^\alpha$ for the set of all injective sequences $t \colon \alpha \to \NN$ that are \emph{co-infinite}, i.e.\ such that $|\NN \setminus \operatorname{im}(t)| = \infty$. 
Then $T = \bigcup_{\alpha < \omega_1} T^\alpha$ is a well-founded tree when ordered by \emph{extension}, i.e.\ $t \leq t'$ if $t = t' \restriction \operatorname{dom}(t)$. 
For a sequence $s \in T^{\alpha+1}$ of successor length, let $\operatorname{last}(s) := s(\alpha) \in \NN$ be the last value of $s$; and $s^\star := s \restriction \alpha \in T^\alpha$ its immediate predecessor.
Put $\Sigma(T) = \bigcup_{\alpha < \omega_1} T^{\alpha +1}$. 
For any $t \in T$, let 
$$A_t := \set{s \leq t}:{s \in \Sigma(T), \; \operatorname{last}(s) = \min \p{ \operatorname{im}(t) \setminus \operatorname{im}(s^\star)}},$$
and let $A^\star_t = \set{s^\star}:{s \in A_t}$.

Let $\mathbf{G}$ be the graph with vertex set $V(\mathbf{G}) = T$ and edge set $E(\mathbf{G}) = \set{t't}:{t' \in A^\star_t}$.

\begin{thm}
The graph $\mathbf{G}$ is uncountably chromatic yet every uncountable set of vertices in $\mathbf{G}$ contains two points which are connected by only finitely many independent paths in $\mathbf{G}$.
\end{thm}

\section{The proof}

We first show that every uncountable set of vertices $A \subseteq V(\mathbf{G})$ contains two points which are connected by only finitely many independent paths in $\mathbf{G}$. For $s \in T$ we write $s \!\downarrow := \set{t \in T}:{t < s}$, and note that the definition of $A_s$ implies that for all $ s \leq u \in T$ we have 
$$A_{u} \cap s\!\downarrow \; \subseteq A_s \; \text{ and } \; A^\star_u \cap s\!\downarrow \; \subseteq A^\star_s. \quad \quad \quad (\star)$$
Since $T$ contains no uncountable chains, the set $A$ contains two vertices $t$ and $t'$ that are incomparable in $T$.
Let $\alpha \in \operatorname{dom}(t)$ be minimal such that $t(\alpha) \neq t'(\alpha)$, and consider $s = t \restriction (\alpha+ 1)$. Then $(\star)$ implies that every $t-t'$ path meets $A^\star_s$, and since $|A^\star_s| = |A_s| \leq \operatorname{last}(s)$ is finite, there are only finitely many independent $t-t'$ paths in $\mathbf{G}$.

\medskip
It remains to show that $\mathbf{G}$ has chromatic number $\chi(\mathbf{G}) =\aleph_1$. Colouring the elements of each $T^\alpha$ with a new colour shows $\chi(\mathbf{G}) \leq \aleph_1$. To see  $\chi(\mathbf{G}) \geq \aleph_1$, suppose for a contradiction that $c\colon V(\mathbf{G}) \to \NN$ is a proper colouring. 

For $t \in \Sigma(T)$, we say $t' \in T$ is an \emph{extension of $t$} if $t' \geq t$ and $\operatorname{im}(t') \setminus \operatorname{im}(t) \subseteq \Set{n \in \NN \colon n > \operatorname{last}(t)}$. We say $t' \in T$ is a \emph{1-extension of $t$} if it has the stronger property that $t' > t$ and, letting $a_1$ be the minimal element of $\NN \setminus \operatorname{im}(t)$ with $a_1 > \operatorname{last}(t)$, we have $\operatorname{im}(t') \setminus \operatorname{im}(t) \subseteq \Set{n \in \NN \colon n > a_{1}}$. In this case we also say that the $1$-extension $t'$ \emph{skips} $a_1$.

\begin{clm} Every $t$ in $\Sigma(T)$ has an extension $t'$ in $\Sigma(T)$ such that every $1$-extension $t''$ of $t'$ satisfies $c(t'') > c(t^\star)$.
\end{clm}

Suppose for a contradiction that the claim is false. Then there exists $t_0 \in \Sigma(T)$ such that all its extensions $t' \in \Sigma(T)$ have a $1$-extension $t''$ such that $c(t'') \leq c(t_0^\star)$.
Then for the extension $t'_0 =t_0$ of $t_0$, there is a 1-extension  $t''_0$ of $t'_0$ skipping $a_1 \in \NN$ with $c(t''_0) \leq c(t^\star_0)$. Let $t'_1 := t''_0 {}^\frown a_1$. Then $t'_1 \in \Sigma(T)$ is itself an extension of $t_0$, so it has a 1-extension $t''_1$ skipping $a_2$ with $c(t''_1) \leq c(t^\star_0)$. Let $t'_2 := t''_1 {}^\frown a_2$. And so on.  

Now $a_{m+1}$ witnesses that $t'_{m+1} \in A_{t''_n}$, and so $t''_{m} \in A^\star_{t''_n}$ whenever $m < n \in \NN$. 
Hence, the vertices $\set{t''_n}:{n \in \NN}$ induce a complete subgraph of $\mathbf{G}$, contradicting that they been coloured using only colours $\leq c(t^\star_0)$. This proves the claim.

\medskip
We now complete the proof as follows: Fix an arbitrary $t_0 \in \Sigma(T)$. 
Let $t'_0 \in \Sigma(T)$ be an extension of $t_0$ as in the claim. 
Let $a_1<a_2$ be the two smallest elements of $\NN \setminus \operatorname{im}(t'_0)$ above $\operatorname{last}(t'_0)$.
Let $t_1 := t'_0 {}^\frown a_2$. 
Let $t'_1$ be an extension of $t_1$  as in the claim.
Let $a_3<a_4$ be the two smallest elements of $\NN \setminus \operatorname{im}(t'_1)$ above $\operatorname{last}(t'_1)$.
Let $t_2:= t'_1 {}^\frown a_4$. 
And so on. 

Then $\hat{t} = \bigcup_{n \in \NN} t_n$ is an injective sequence.  
Moreover, $a_1,a_3,a_5,\ldots$ witness that $\hat{t}$ is co-infinite, giving $\hat{t} \in T$. 
But for each $n \in \NN$, the sequence $\hat{t}$ is a 1-extension (skipping $a_{2n+1}$) 
of the extension $t'_n$ of $t_n$, so 
$c(t^\star_n) \leq c(\hat{t})$
according to the claim.  However, $a_2,a_4,a_6,\ldots$ witness that the vertices  $\set{t^\star_n}:{n \in \NN}$ induce a complete subgraph of $\mathbf{G}$,  a contradiction. \hfill \qed

\section{Remarks}

(1) In the terminology of \cite{kurkofka2021representation}, the graph $\mathbf{G}$ is a \emph{$T$-graph of finite adhesion}. The construction of the sets $A_t$ is inspired by an argument from \cite{diestel2001normal}.

(2) The graph with vertex set $T$ but edge set $\set{t't}:{t' < t, \; t' \in A_t}$ has countable chromatic number by colouring all $s \in \Sigma(T)$ by colour $\operatorname{last}(s)$, and noticing that $A_t \subset \Sigma(T)$ for all $t \in T$ implies that $T \setminus \Sigma(T)$ is independent.

(3) The following version of the Erd\H{o}s-Hajnal problem remains open:
\emph{Does every uncountably chromatic graph have a countably infinite, infinitely connected subgraph?}

\bibliographystyle{plain}
\bibliography{reference}

\end{document}